\documentclass[reqno]{amsart}

\numberwithin{Theorem}{section}
\numberwithin{equation}{section}

\allowdisplaybreaks

\begin{document}

\title[A simple proof of the $_6\psi_6$ sum]
{A simple proof of Bailey's very-well-poised 
$\hbox{${}_{\boldsymbol 6}\boldsymbol\psi_{\boldsymbol 6}$}$ summation}
\author{Michael Schlosser}

\address{Department of Mathematics, The Ohio State University,
231 West 18th Avenue, Columbus, Ohio 43210, USA}
\email{mschloss@math.ohio-state.edu}
\date{September 25, 2000}


\begin{abstract}
We give elementary derivations of some classical
summation formulae for bilateral (basic)
hypergeometric series. In particular, we
apply Gau\ss{}' $_2F_1$ summation and elementary
series manipulations to give a simple proof of
Dougall's $_2H_2$ summation. Similarly, we apply
Rogers' nonterminating $_6\phi_5$ summation
and elementary series manipulations to give a simple proof of Bailey's
very-well-poised $_6\psi_6$ summation.
Our method of proof extends M.~Jackson's first elementary
proof of Ramanujan's $_1\psi_1$ summation.
\end{abstract}

\maketitle

\section{Introduction}

The theories of unilateral
(or one-sided) hypergeometric and basic ($q$-)hyper\-geometric series
have quite a rich history dating back to at least Euler.
Formulae for {\em bilateral} (basic) hypergeometric
series were not discovered until 1907 when Dougall~\cite{dougall},
using residue calculus, derived summations for the bilateral $_2H_2$
and very-well-poised $_5H_5$ series.
Ramanujan~\cite{hardy} extended the $q$-binomial theorem by finding
a summation formula for the bilateral $_1\psi_1$ series.
Later, Bailey~\cite{bail66},\cite{bail22} carried out
systematical investigations of summations and transformations
for bilateral basic hypergeometric series.
Further significant contributions were made by
Slater~\cite{slatertf},\cite{slater}, a student of Bailey.
See \cite{grhyp} and \cite{slater} for an excellent survey of the
above classical material.

Bailey's~\cite[Eq.~(4.7)]{bail66} very-well-poised $_6\psi_6$ summation
(cf.~\cite[Eq.~(5.3.1)]{grhyp}) is a very powerful identity,
as it stands at the top of the classical
hierarchy of summation formulae for bilateral series.
Some of the applications of the $_6\psi_6$ summation to partitions and
number theory are given in Andrews~\cite{andappl}.
Though several proofs of Bailey's $_6\psi_6$ summation are already
known (see, e.~g., Bailey~\cite{ bail66}, Slater and Lakin~\cite{slatlat},
Andrews~\cite{andappl}, Askey and Ismail~\cite{askmail}, and
Askey~\cite{askeyII}), none of them is entirely elementary.
Here we provide a new simple proof of the very-well-poised
$_6\psi_6$ summation formula, directly from three applications of
Rogers'~\cite[p.~29, second eq.]{rogers} nonterminating
$_6\phi_5$ summation (cf.~\cite[Eq.~(2.7.1)]{grhyp}) and elementary
manipulations of series.

The method of proof we apply extends that already used by
M.~Jackson~\cite[Sec.~4]{mjack} in her first elementary proof
(as pointed out to us by George Andrews~\cite{andpriv}) of Ramanujan's
$_1\psi_1$ summation formula~\cite{hardy} (cf.~\cite[Eq.~(5.2.1)]{grhyp}).
Jackson's proof essentially derives the $_1\psi_1$ summation from
the $q$-Gau\ss{} summation, by manipulation of series.
In view of this background,
it is surprising that this method has not been
further applied for half a century.
A possible explanation is that the applicability of her method was viewed
as too limited. In fact, only after changing the order of steps
in Jackson's proof, we were able to extend her proof to a ``method". 

Indeed, the method can also be applied
to derive other summations.
After recalling some notation for (basic) hypergeometric series
in Section~\ref{sec0}, we review Jackson's elementary proof of
the $_1\psi_1$ summation in Section~\ref{sec1}.
In Section~\ref{sec2}, we apply our extension of Jackson's method to
give an elementary proof of Dougall's~\cite{dougall}
$_2H_2$ summation.
Finally, in Section~\ref{sec6}, we give an elementary derivation of
Bailey's very-well-poised $_6\psi_6$ summation.

We want to point out that by using a similar but slightly different method,
the author~\cite{schleltf} has found elementary
derivations of transformations for bilateral basic hypergeometric series.
In fact, in \cite{schleltf} we use Bailey's~\cite{bail66}
nonterminating very-well-poised $_8\phi_7$ summation theorem
combined with bilateral series identities to derive
a very-well-poised $_8\psi_8$ transformation, a
very-well-poised $_{10}\psi_{10}$ transformation, and 
by induction, Slater's~\cite{slatertf} general transformation
for very-well-poised $_{2r}\psi_{2r}$ series.
Similarly, some other bilateral series identities are also
elementarily derived in \cite{schleltf}.

In the near future, we plan to apply the methods of this article and
of \cite{schleltf} to the settings of multiple basic hypergeometric
series. See Milne~\cite{milne}, Gustafson~\cite{gus},
v.~Diejen~\cite{diejen}, and Schlosser~\cite{hypdet},
for several of these different settings. We are quite confident that
we may not only get simpler proofs for already known results but should
also obtain derivations of new formulae.

Finally, we wish to gratefully acknowledge the helpful comments and
suggestions of George Andrews, Mourad Ismail, and Stephen Milne.

\section{Background and notation}\label{sec0}

Here we recall some notation for hypergeometric series (cf.~\cite{slater}),
and basic hypergeometric series (cf.~\cite{grhyp}).

We define the {\em shifted factorial} for all integers $k$ by
the following quotient of Gamma functions (cf.~\cite[Sec.~1.1]{sfaar}),
\begin{equation*}
(a)_k:=\frac{\Gamma(a+k)}{\Gamma(a)}.
\end{equation*}
Further, the (ordinary) {\em hypergeometric $_rF_s$ series} is defined as
\begin{equation}\label{defohyp}
{}_rF_s\!\left[\begin{matrix}a_1,a_2,\dots,a_r\\
b_1,b_2,\dots,b_s\end{matrix};z\right]:=
\sum _{k=0} ^{\infty}\frac {(a_1)_k\dots(a_r)_k}
{(b_1)_k\dots(b_s)_k}\frac{z^k}{k!},
\end{equation}
and the {\em bilateral hypergeometric $_rH_s$ series} as
\begin{equation}\label{defohypb}
{}_rH_s\!\left[\begin{matrix}a_1,a_2,\dots,a_r\\
b_1,b_2,\dots,b_s\end{matrix};z\right]:=
\sum _{k=-\infty} ^{\infty}\frac {(a_1)_k\dots(a_r)_k}
{(b_1)_k\dots (b_s)_k}z^k.
\end{equation}
See \cite[p.~45 and p.~181]{slater} for the criteria
of when these series terminate, or, if not, when they converge. 

Let $q$ be a complex number such that $0<|q|<1$. We define the
{\em $q$-shifted factorial} for all integers $k$ by
\begin{equation*}
(a;q)_{\infty}:=\prod_{j=0}^{\infty}(1-aq^j)\qquad\text{and}\qquad
(a;q)_k:=\frac{(a;q)_{\infty}}{(aq^k;q)_{\infty}}.
\end{equation*}
For brevity, we employ the usual notation
\begin{equation*}
(a_1,\ldots,a_m;q)_k\equiv (a_1;q)_k\dots(a_m;q)_k
\end{equation*}
where $k$ is an integer or infinity. Further, we utilize the notations
\begin{equation}\label{defhyp}
_r\phi_s\!\left[\begin{matrix}a_1,a_2,\dots,a_r\\
b_1,b_2,\dots,b_s\end{matrix};q,z\right]:=
\sum _{k=0} ^{\infty}\frac {(a_1,a_2,\dots,a_r;q)_k}
{(q,b_1,\dots,b_s;q)_k}\left((-1)^kq^{\binom k2}\right)^{1+s-r}z^k,
\end{equation}
and
\begin{equation}\label{defhypb}
_r\psi_s\!\left[\begin{matrix}a_1,a_2,\dots,a_r\\
b_1,b_2,\dots,b_s\end{matrix};q,z\right]:=
\sum _{k=-\infty} ^{\infty}\frac {(a_1,a_2,\dots,a_r;q)_k}
{(b_1,b_2,\dots,b_s;q)_k}\left((-1)^kq^{\binom k2}\right)^{s-r}z^k,
\end{equation}
for {\em basic hypergeometric $_r\phi_s$ series}, and {\em bilateral basic
hypergeometric $_r\psi_s$
series}, respectively. See \cite[p.~25 and p.~125]{grhyp} for the criteria
of when these series terminate, or, if not, when they converge. 

We want to point out that many theorems for ${}_rF_s$ or ${}_rH_s$
series can be obtained by considering certain ``$q\to 1$ limiting cases"
of corresponding theorems for ${}_r\phi_s$ or ${}_r\psi_s$ series,
respectively. For instance, we describe such a $q\to 1$ limiting case
after stating the $q$-binomial theorem in \eqref{10gl}. 
A similar $q\to 1$ limiting case leads from the $q$-Gau\ss{}
summation~\eqref{21gl} to the ordinary Gau\ss{}
summation~\eqref{2f1gl}. The situation is different
for the $_{1}\psi_1$ series, though. We have a summation for the general
$_{1}\psi_1$, but not for the $_{1}H_1$.
On the other hand, the general $_2H_2$ with unit argument
is summable but the general
$_2\psi_2$ is not. Many theorems for very-well-poised
$_{r+1}\phi_r$ series can be specialized to theorems
for very-well-poised $_rF_{r-1}$ series. For the notion of ({\em
very-}){\em well-poised}, see \cite[Sec.~2.1]{grhyp}. 
For detailed treatises on hypergeometric and basic hypergeometric
series, we refer to Slater~\cite{slater}, and Gasper and Rahman~\cite{grhyp}.

In our computations in the following sections, we make heavily use of some
elementary identities involving ($q$-)shifted factorials which are listed
in Slater~\cite[Appendix~I]{slater}, and
Gasper and Rahman~\cite[Appendix~I]{grhyp}.

\section{M.~Jackson's proof of Ramanujan's $_1\psi_1$ summation}\label{sec1}

The $q$-binomial theorem,
\begin{equation}\label{10gl}
{}_1\phi_0\!\left[\begin{matrix}a\\
-\end{matrix};q,z\right]=
\frac{(az;q)_{\infty}}{(z;q)_{\infty}},
\end{equation}
where the series either terminates, or $|z|<1$, for convergence,
was first discovered by Cauchy~\cite{cauchy} (cf.~\cite[Sec.~1.3]{grhyp}).
It reduces to the ordinary binomial theorem as $a\mapsto q^a$ and $q\to1^-$.

A bilateral extension of the $q$-binomial theorem \eqref{10gl},
the $_1\psi_1$ summation, was found by the legendary Indian
mathematician Ramanujan~\cite{hardy} (cf.~\cite[Eq.~(5.2.1)]{grhyp}).
It reads as follows:
\begin{equation}\label{11gl}
{}_1\psi_1\!\left[\begin{matrix}a\\
b\end{matrix};q,z\right]=
\frac{(q,b/a,az,q/az;q)_{\infty}}{(b,q/a,z,b/az;q)_{\infty}},
\end{equation}
where the series either terminates, or $|b/a|<|z|<1$, for convergence.
Clearly, \eqref{11gl} reduces to \eqref{10gl} when $b=q$.

Unfortunately, Ramanujan did not provide a proof for his summation
formula. Hahn~\cite[$\kappa=0$ in Eq.~(4.7)]{hahn} independently
established \eqref{11gl} by considering a first order homogeneous
$q$-difference equation. Hahn thus published the first proof
of the $_1\psi_1$ summation. Not much later, M.~Jackson~\cite[Sec.~4]{mjack}
gave the first elementary proof of \eqref{11gl}.
Her proof derives the $_1\psi_1$ summation from the $q$-Gau\ss{}
summation, by manipulation of series. It turns out that
Jackson's method is effective for proving also other bilateral
summation formulae. Since Jackson's short proof of
Ramanujan's $_1\psi_1$  summation seems to be not so well known,
we review her proof in the following.

Before we continue, we want to point out that
there are also many other nice proofs of the $_1\psi_1$ summation
in the literature. A simple and elegant
proof of the $_1\psi_1$ summation formula was given by Ismail~\cite{ismail}
who showed that the $_1\psi_1$ summation is an immediate consequence of
the $q$-binomial theorem and analytic continuation.

M.~Jackson's elementary proof of \eqref{11gl}
makes use of a suitable specialization of
Heine's~\cite{heine} $q$-Gau\ss{} summation (cf.~\cite[Eq.~(II.8)]{grhyp}),
\begin{equation}\label{21gl}
{}_2\phi_1\!\left[\begin{matrix}a,b\\
c\end{matrix};q,\frac c{ab}\right]=
\frac{(c/a,c/b;q)_{\infty}}{(c,c/ab;q)_{\infty}},
\end{equation}
where the series either terminates, or $|c/ab|<1$, for convergence.

In \eqref{21gl}, we perform the substitutions $a\mapsto aq^n$,
$b\mapsto q/b$, and $c\mapsto q^{1+n}$, and obtain
\begin{equation}\label{21gl1}
{}_2\phi_1\!\left[\begin{matrix}aq^n,q/b\\
q^{1+n}\end{matrix};q,\frac ba\right]=
\frac{(q/a,bq^n;q)_{\infty}}{(q^{1+n},b/a;q)_{\infty}},
\end{equation}
provided $|b/a|<1$.

Using some elementary identities for $q$-shifted factorials
(see, e.~g., Gasper and Rahman~\cite[Appendix~I]{grhyp})
we can rewrite equation~\eqref{21gl1} as
\begin{equation}\label{key0}
\frac {(q,b/a;q)_{\infty}}
{(q/a,b;q)_{\infty}}
\sum_{k=0} ^{\infty}\frac {(q/b;q)_k(a;q)_{n+k}}
{(q;q)_k(q;q)_{n+k}}
\left(\frac {b} {a}\right)^k
=\frac {(a;q)_n}{(b;q)_n}.
\end{equation}
In this identity, we multiply both sides by $z^n$
and sum over all integers $n$.

On the right side we obtain
\begin{equation*}
{}_1\psi_1\!\left[\begin{matrix}a\\
b\end{matrix};q,z\right].
\end{equation*}
On the left side we obtain
\begin{equation}\label{21gl2}
\frac {(q,b/a;q)_{\infty}}
{(q/a,b;q)_{\infty}}\sum_{n=-\infty}^{\infty}z^n
\sum_{k=0} ^{\infty}\frac {(q/b;q)_k(a;q)_{n+k}}
{(q;q)_k(q;q)_{n+k}}
\left(\frac {b} {a}\right)^k.
\end{equation}
Next, we interchange summations in \eqref{21gl2} and shift
the inner index $n\mapsto n-k$.
(Observe that the sum over $n$ is terminated by the term
$(q;q)_{n+k}^{-1}$ from below.)  We obtain
\begin{equation*}
\frac {(q,b/a;q)_{\infty}}
{(q/a,b;q)_{\infty}}
\sum_{k=0} ^{\infty}\frac {(q/b;q)_k}{(q;q)_k}
\left(\frac {b} {az}\right)^k
\sum_{n=0}^{\infty}\frac {(a;q)_n}{(q;q)_{n}}z^n.
\end{equation*}
Now, twice application of the $q$-binomial theorem \eqref{10gl}
gives us the right side of \eqref{11gl}, as desired.

Now, we have to admit that M.~Jackson did not give her proof
in the above precise order. In fact, her proof
in~\cite[Sec.~4]{mjack} goes backwards. (This is also how the author
originally rediscovered Jackson's proof.) She started with
the $_1\psi_1$ summation~\eqref{11gl} and equated coefficients
of $z^n$ on both sides. The resulting identity is true by the
$q$-Gau\ss{} summation.

A reason why M.~Jackson's method of proof has so far not been used
to prove other
bilateral summations could be that the applicability of
her derivation was viewed as too limited. Equating
coefficients of a power of a Laurent series variable in a
bilateral basic hypergeometric series identity
is easy if, as in \eqref{11gl}, there is an argument $z$ which is
independent of the other parameters. But this seems to be
more particluar to the $_1\psi_1$ summation, as not in all
bilateral series there is such an independent argument.
The starting point for making M.~Jackson's proof to a ``method"
is to read the proof backwards, as displayed above.
The essence here is that a unilateral series identity, \eqref{21gl},
is specialized such that there is the factor $(q;q)_{n+k}^{-1}$
in the series, see \eqref{key0},
so that summing over all $n$ again gives a (summable) unilateral series. 

In the next two sections, we use the method to give proofs of two
other important bilateral hypergeometric and basic
hypergeometric summation theorems.
In particular, in Section~\ref{sec2}, we give a simple proof of 
Dougall's $_2H_2$ summation, whereas in Section~\ref{sec6},
we give a simple proof of Bailey's very-well-poised $_6\psi_6$ summation.

\section{Dougall's $_2H_2$ summation}\label{sec2}

In Section~\ref{sec1}, we multiplied both sides of the identity~\eqref{key0}
by a suitable factor depending on $n$ and summed over all integers $n$.
On one side, we interchanged sums and found that the inner sum was
summable by the $q$-binomial theorem. Now, what if we start with
a different factor following a similar procedure
such that we can evaluate the inner sum by, say,
the $q$-Gau\ss{} summation? If the analysis works out we may end up
with an evaluation for a $_2\psi_2$ series. Let us see what happens:

In identity~\eqref{key0}, let us first replace $b$ by $c$.
Then  we multiply both sides by
\begin{equation*}
\frac{(b;q)_n}{(d;q)_n}\left(\frac{d}{ab}\right)^n
\end{equation*}
and sum over all integers $n$.

On the right side we obtain
\begin{equation}\label{22tf1}
{}_2\psi_2\!\left[\begin{matrix}a,b\\
c,d\end{matrix};q,\frac{d}{ab}\right].
\end{equation}
On the left side we obtain
\begin{equation}\label{22gl1}
\frac {(q,c/a;q)_{\infty}}
{(q/a,c;q)_{\infty}}\sum_{n=-\infty}^{\infty}\frac{(b;q)_n}{(d;q)_n}
\left(\frac{d}{ab}\right)^n
\sum_{k=0} ^{\infty}\frac {(q/c;q)_k(a;q)_{n+k}}
{(q;q)_k(q;q)_{n+k}}
\left(\frac {c} {a}\right)^k.
\end{equation}
Next, we interchange summations in \eqref{22gl1} and shift
the inner index $n\mapsto n-k$.  We obtain, again using some elementary
identities for $q$-shifted factorials,
\begin{equation*}
\frac {(q,c/a;q)_{\infty}}
{(q/a,c;q)_{\infty}}
\sum_{k=0} ^{\infty}\frac {(q/c;q)_k(b;q)_{-k}}{(q;q)_k(d;q)_{-k}}
\left(\frac {bc} {d}\right)^k
\sum_{n=0}^{\infty}\frac {(a,bq^{-k};q)_n}{(q,dq^{-k};q)_{n}}
\left(\frac{d}{ab}\right)^n.
\end{equation*}
Now the inner sum, provided $|d/ab|<1$, can
be evaluated by~\eqref{21gl} and we obtain
\begin{equation*}
\frac {(q,c/a;q)_{\infty}}
{(q/a,c;q)_{\infty}}
\sum_{k=0} ^{\infty}\frac {(q/c;q)_k(b;q)_{-k}}{(q;q)_k(d;q)_{-k}}
\left(\frac {bc} {d}\right)^k
\frac{(dq^{-k}/a,d/b;q)_{\infty}}{(dq^{-k},d/ab;q)_{\infty}},
\end{equation*}
which can be simplied to
\begin{equation}\label{22tf2}
\frac {(q,c/a,d/a,d/b;q)_{\infty}}
{(q/a,c,d,d/ab;q)_{\infty}}
\sum_{k=0} ^{\infty}\frac {(q/c,aq/d;q)_k}{(q,q/b;q)_k}
\left(\frac {c} {a}\right)^k.
\end{equation}
Hence, equating \eqref{22tf1} and \eqref{22tf2},
we have derived the transformation
\begin{equation}\label{22tf}
{}_2\psi_2\!\left[\begin{matrix}a,b\\
c,d\end{matrix};q,\frac{d}{ab}\right]=
\frac {(q,c/a,d/a,d/b;q)_{\infty}}
{(q/a,c,d,d/ab;q)_{\infty}}\;
{}_2\phi_1\!\left[\begin{matrix}q/c,aq/d\\
q/b\end{matrix};q,\frac{c}a\right],
\end{equation}
where the series terminate, or $\max(|d/ab|,|c|,|c/a|)<1$, for convergence.
Unfortunately, the $_2\phi_1$ on the right side of \eqref{22tf}
simplifies only in special cases. If $d=aq$, then the $_2\phi_1$ sum
reduces just to the first term, 1, and we have the summation
\begin{equation}\label{22sum}
{}_2\psi_2\!\left[\begin{matrix}a,b\\
aq,c\end{matrix};q,\frac{q}{b}\right]=
\frac {(q,q,aq/b,c/a;q)_{\infty}}
{(aq,q/a,q/b,c;q)_{\infty}},
\end{equation}
where the series terminates, or $\max(|q/b|,|c|)<1$, for convergence.

We want to add that the transformation in \eqref{22tf} is a special case
of Bailey's~\cite[Eq.~(2.3)]{bail22} $_2\psi_2$ transformation,
\begin{equation}\label{22btf}
{}_2\psi_2\!\left[\begin{matrix}a,b\\
c,d\end{matrix};q,z\right]=
\frac {(az,d/a,c/b,dq/abz;q)_{\infty}}
{(z,d,q/b,cd/abz;q)_{\infty}}\;
{}_2\psi_2\!\left[\begin{matrix}a,abz/d\\
az,c\end{matrix};q,\frac da\right],
\end{equation}
where the series terminate, or $\max(|z|,|cd/abz|,|d/a|,|c/b|)<1$,
for convergence.
Namely, if we perform in \eqref{22btf} the simultaneous substitutions
$a\mapsto b$, $b\mapsto a$, and $z\mapsto d/ab$,
and reverse the order of summation in the truncated series
on the right side, we obtain \eqref{22tf}.

In Section~\ref{sec1}, we found, following M.~Jackson, a sum
for a general $_1\psi_1$ series. So far in this section, we applied her
method to obtain a transformation for a particular $_2\psi_2$ into a
(multiple of a) $_2\phi_1$ series. As a matter of fact,
there is no closed form (as a product of linear factors) for the
summation of a general $_2\psi_2$ series.
The situation is different in the $q\to 1$ case, though.

In the following, we review the classical $_2F_1$ and $_2H_2$
summations and then prove the latter by our elementary method.

In his doctoral dissertation~\cite{gauss}, Gau\ss{} showed that
\begin{equation}\label{2f1gl}
{}_2F_1\!\left[\begin{matrix}a,b\\
c\end{matrix};1\right]=
\frac{\Gamma(c)\Gamma(c-a-b)}{\Gamma(c-a)\Gamma(c-b)},
\end{equation}
where the series either terminates, or $\Re(c-a-b)>0$, for convergence.

Dougall~\cite[Sec.~13]{dougall} extended this result to
\begin{equation}\label{2h2gl}
{}_2H_2\!\left[\begin{matrix}a,b\\
c,d\end{matrix};1\right]=
\frac{\Gamma(1-a)\Gamma(1-b)\Gamma(c)\Gamma(d)\Gamma(c+d-a-b-1)}
{\Gamma(c-a)\Gamma(c-b)\Gamma(d-a)\Gamma(d-b)},
\end{equation}
where the series either terminates, or $\Re(c+d-a-b-1)>0$, for convergence.
Clearly, the $d\to 1$ case of \eqref{2h2gl} is \eqref{2f1gl}.

We are ready to derive \eqref{2h2gl} from \eqref{2f1gl}:
In \eqref{2f1gl}, we perform the simultaneous substitutions
$a\mapsto a+n$, $b\mapsto 1-c$, and $c\mapsto 1+n$, and obtain
\begin{equation}\label{2f1gl1}
{}_2F_1\!\left[\begin{matrix}a+n,1-c\\
1+n\end{matrix};1\right]=
\frac{\Gamma(1+n)\Gamma(c-a)}{\Gamma(1-a)\Gamma(c+n)},
\end{equation}
provided $\Re(c-a)>0$.

Using some elementary identities for shifted factorials (see, e.~g.,
Slater~\cite[Appendix~I]{slater}) we can rewrite equation~\eqref{2f1gl1} as
\begin{equation}\label{2f1gl2}
\frac{\Gamma(1-a)\Gamma(c)}{\Gamma(c-a)}
\sum_{k=0}^{\infty}\frac{(1-c)_k(a)_{n+k}}{(1)_k(1)_{n+k}}
=\frac{(a)_n}{(c)_n}.
\end{equation}
Alternatively, we could have used \eqref{key0} with the substitutions
$a\mapsto q^{a}$ and $b\mapsto q^{c}$, and then let $q\to 1^-$,
to arrive directly at \eqref{2f1gl2}.

In \eqref{2f1gl2}, we multiply both sides by $(b)_n/(d)_n$
and sum over all integers $n$.

On the right side we obtain
\begin{equation*}
{}_2H_2\!\left[\begin{matrix}a,b\\
c,d\end{matrix};1\right].
\end{equation*}
On the left side we obtain
\begin{equation}\label{2f1gl3}
\frac{\Gamma(1-a)\Gamma(c)}{\Gamma(c-a)}
\sum_{n=-\infty}^{\infty}\frac{(b)_n}{(d)_n}
\sum_{k=0}^{\infty}\frac{(1-c)_k(a)_{n+k}}{(1)_k(1)_{n+k}}.
\end{equation}
Next, we interchange summations in \eqref{2f1gl3} and shift
the inner index $n\mapsto n-k$.
(Observe that the sum over $n$ is terminated by the term
$(1)_{n+k}^{-1}$ from below.)  We obtain
\begin{equation*}
\frac{\Gamma(1-a)\Gamma(c)}{\Gamma(c-a)}
\sum_{k=0}^{\infty}\frac{(1-c)_k(b)_{-k}}{(1)_k(d)_{-k}}
\sum_{n=0}^{\infty}\frac{(a)_n(b-k)_n}{(1)_n(d-k)_n}.
\end{equation*}
Now, the inner sum, provided $\Re(d-a-b)>0$, can be evaluated by
\eqref{2f1gl} and we obtain
\begin{equation*}
\frac{\Gamma(1-a)\Gamma(c)}{\Gamma(c-a)}
\sum_{k=0}^{\infty}\frac{(1-c)_k(b)_{-k}}{(1)_k(d)_{-k}}
\frac{\Gamma(d-k)\Gamma(d-a-b)}{\Gamma(d-a-k)\Gamma(d-b)},
\end{equation*}
which can be simplified to
\begin{equation*}
\frac{\Gamma(1-a)\Gamma(c)\Gamma(d)\Gamma(d-a-b)}
{\Gamma(c-a)\Gamma(d-a)\Gamma(d-b)}
\sum_{k=0}^{\infty}\frac{(1-c)_k(1+a-d)_k}{(1)_k(1-b)_k}.
\end{equation*}
To the last inner sum, provided $\Re(c+d-a-b-1)>0$,
we can again apply \eqref{2f1gl} and 
eventually obtain the right side of \eqref{2h2gl}, as desired.

We note here that to apply Gau\ss{}' $_2F_1$ summation theorem
three times, we
needed certain conditions of the parameters, for convergence.
These were $\Re(c-a)>0$, $\Re(d-a-b)>0$, and $\Re(c+d-a-b-1)>0$.
But in the end the first two of these conditions may be removed
by analytic continuation. 
In particular, both sides of identity \eqref{2h2gl} are analytic in
$a$ for $\Re(a)<\Re(c+d-b-1)$ (and excluding some poles).
In the course of our derivation, we have shown
the identity for $\Re(a)<\min(\Re(c),\Re(d-b),\Re(c+d-b-1))$
(excluding some poles).
By analytic continuation, we extend the identity, when defined,
to be valid for $\Re(a)<\Re(c+d-b-1)$, the region of convergence of
the series.

\section{Bailey's very-well-poised $_6\psi_6$ summation}\label{sec6}

One of the most powerful identities for bilateral basic hypergeometric series
is Bailey's very-well-poised $_6\psi_6$ summation:
\begin{multline}\label{66gl}
{}_6\psi_6\!\left[\begin{matrix}q\sqrt{a},-q\sqrt{a},b,c,d,e\\
\sqrt{a},-\sqrt{a},aq/b,aq/c,aq/d,aq/e\end{matrix};q,
\frac{a^2q}{bcde}\right]\\
=\frac {(aq,aq/bc,aq/bd,aq/be,aq/cd,aq/ce,aq/de,q,q/a;q)_{\infty}}
{(aq/b,aq/c,aq/d,aq/e,q/b,q/c,q/d,q/e,a^2q/bcde;q)_{\infty}},
\end{multline}
provided the series either terminates, or $|a^2q/bcde|<1$,
for convergence.

To prove Bailey's $_6\psi_6$ summation, we start with a suitable
specialization of Rogers' $_6\phi_5$
summation:
\begin{equation}\label{65gl}
{}_6\phi_5\!\left[\begin{matrix}a,\,q\sqrt{a},-q\sqrt{a},b,c,d\\
\sqrt{a},-\sqrt{a},aq/b,aq/c,aq/d\end{matrix};q,
\frac{aq}{bcd}\right]
=\frac {(aq,aq/bc,aq/bd,aq/cd;q)_{\infty}}
{(aq/b,aq/c,aq/d,aq/bcd;q)_{\infty}},
\end{equation}
provided the series either terminates, or $|aq/bcd|<1$, for
convergence. Note that \eqref{65gl} is just the special case $e\mapsto a$
of \eqref{66gl}.

In \eqref{65gl}, we perform the simultaneous substitutions $a\mapsto c/a$,
$b\mapsto b/a$, $c\mapsto cq^{n}$ and $d\mapsto cq^{-n}/a$, and obtain
\begin{multline}\label{65gl1}
{}_6\phi_5\!\left[\begin{matrix}c/a,\,q\sqrt{c/a},-q\sqrt{c/a},b/a,cq^n,
cq^{-n}/a\\
\sqrt{c/a},-\sqrt{c/a},cq/b,q^{1-n}/a,q^{1+n}\end{matrix};q,
\frac{aq}{bc}\right]\\
=\frac {(cq/a,q^{1-n}/b,aq^{1+n}/b,q/c;q)_{\infty}}
{(cq/b,q^{1-n}/a,q^{1+n},aq/bc;q)_{\infty}},
\end{multline}
where $|aq/bc|<1$.

Using some elementary identities for $q$-shifted factorials
(see, e.~g., Gasper and Rahman~\cite[Appendix~I]{grhyp})
we can rewrite equation~\eqref{65gl1} as
\begin{multline}\label{key1}
\frac {(cq/b,q/a,q,aq/bc;q)_{\infty}}
{(cq/a,q/b,aq/b,q/c;q)_{\infty}}
\sum_{k=0} ^{\infty}\frac {(1-cq^{2k}/a)} {(1-c/a)}
\frac {(c/a,b/a;q)_k(c;q)_{n+k}(a;q)_{n-k}}
{(q,cq/b;q)_k(q;q)_{n+k}(aq/c;q)_{n-k}}
\left(\frac {a} {b}\right)^k\\
=\frac {(b,c;q)_n}{(aq/b,aq/c;q)_n}\left(\frac {a} {b}\right)^n.
\end{multline}
In this identity, we multiply both sides by
\begin{equation*}
\frac{(1-aq^{2n})}{(1-a)}\frac{(d,e;q)_n}{(aq/d,aq/e;q)_n}
\left(\frac {aq} {cde}\right)^n
\end{equation*}
and sum over all integers $n$.

On the right side we obtain
\begin{equation*}
{}_6\psi_6\!\left[\begin{matrix}q\sqrt{a},-q\sqrt{a},b,c,d,e\\
\sqrt{a},-\sqrt{a},aq/b,aq/c,aq/d,aq/e\end{matrix};q,
\frac{a^2q}{bcde}\right].
\end{equation*}
On the left side we obtain
\begin{multline}\label{65gl2}
\frac {(cq/b,q/a,q,aq/bc;q)_{\infty}}
{(cq/a,q/b,aq/b,q/c;q)_{\infty}}
\sum_{n=-\infty} ^{\infty}\frac {(1-aq^{2n})} {(1-a)}
\frac {(d,e;q)_n} {(aq/d,aq/e;q)_n}
\left(\frac {aq} {cde}\right)^n\\
\times\sum_{k=0} ^{\infty}\frac {(1-cq^{2k}/a)} {(1-c/a)}
\frac {(c/a,b/a;q)_k(c;q)_{n+k}(a;q)_{n-k}}
{(q,cq/b;q)_k(q;q)_{n+k}(aq/c;q)_{n-k}}
\left(\frac {a} {b}\right)^k.
\end{multline}
Next, we interchange summations in \eqref{65gl2} and shift
the inner index $n\mapsto n-k$.
(Observe that the sum over $n$ is terminated by the term
$(q;q)_{n+k}^{-1}$ from below.)  We obtain, again using some elementary
identities for $q$-shifted factorials,
\begin{multline*}
\frac {(cq/b,q/a,q,aq/bc;q)_{\infty}}
{(cq/a,q/b,aq/b,q/c;q)_{\infty}}
\sum_{k=0} ^{\infty}\frac {(1-cq^{2k}/a)} {(1-c/a)}
\frac {(c/a,b/a;q)_k}{(q,cq/b;q)_k}\\
\times\frac{(1-aq^{-2k})}{(1-a)}
\frac{(a;q)_{-2k}(d,e;q)_{-k}}{(aq/c;q)_{-2k}(aq/d,aq/e;q)_{-k}}
\left(\frac {cde} {bq}\right)^k\\
\times\sum_{n=0} ^{\infty}\frac {(1-aq^{-2k+2n})} {(1-aq^{-2k})}
\frac {(aq^{-2k},c,dq^{-k},eq^{-k};q)_n}
{(q,aq^{1-2k}/c,aq^{1-k}/d,aq^{1-k}/e;q)_n}
\left(\frac {aq} {cde}\right)^n.
\end{multline*}
Now the inner sum, provided $|aq/cde|<1$,
can be evaluated by \eqref{65gl} and we obtain
\begin{multline*}
\frac {(cq/b,q/a,q,aq/bc;q)_{\infty}}
{(cq/a,q/b,aq/b,q/c;q)_{\infty}}
\sum_{k=0} ^{\infty}\frac {(1-cq^{2k}/a)} {(1-c/a)}
\frac {(c/a,b/a;q)_k(aq;q)_{-2k}}
{(q,cq/b;q)_k(aq/c;q)_{-2k}}\\
\times\frac{(d,e;q)_{-k}}{(aq/d,aq/e;q)_{-k}}
\left(\frac {cde} {bq}\right)^k
\frac{(aq^{1-2k},aq^{1-k}/cd,aq^{1-k}/ce,aq/de;q)_{\infty}}
{(aq^{1-2k}/c,aq^{1-k}/d,aq^{1-k}/e,aq/cde;q)_{\infty}},
\end{multline*}
which can be simplified to
\begin{multline*}
\frac {(cq/b,q/a,q,aq/bc,aq,aq/cd,aq/ce,aq/de;q)_{\infty}}
{(cq/a,q/b,aq/b,q/c,aq/c,aq/d,aq/e,aq/cde;q)_{\infty}}\\
\times\sum_{k=0} ^{\infty}\frac {(1-cq^{2k}/a)} {(1-c/a)}
\frac {(c/a,b/a,cd/a,ce/a;q)_k}
{(q,cq/b,q/d,q/e;q)_k}
\left(\frac {a^2q} {bcde}\right)^k.
\end{multline*}
To the last sum, provided $|a^2q/bcde|<1$, we can again apply
\eqref{65gl} and after some simplifications we finally obtain
the right side of \eqref{66gl}, as desired.

Our derivation of the $_6\psi_6$ summation~\eqref{66gl} is simple once
the nonterminating $_6\phi_5$ summation~\eqref{65gl} is given.
But the latter summation follows by an elementary computation from
F.~H.~Jackson's~\cite{jacksum} terminating
$_8\phi_7$ summation (cf.~\cite[Eq.~(2.6.2)]{grhyp})
\begin{multline}\label{87gl}
_8\phi_7\!\left[\begin{matrix}a,\,q\sqrt{a},-q\sqrt{a},b,c,d,
a^2q^{1+n}/bcd,q^{-n}\\
\sqrt{a},-\sqrt{a},aq/b,aq/c,aq/d,bcdq^{-n}/a,aq^{1+n}\end{matrix};q,
q\right]\\
=\frac {(aq,aq/bc,aq/bd,aq/cd;q)_n}
{(aq/b,aq/c,aq/d,aq/bcd;q)_n}
\end{multline}
as $n\to\infty$. Jackson's terminating $_8\phi_7$ summation itself can be
proved by various ways. An algorithmic approach uses the $q$-Zeilberger
algorithm, see Koornwinder~\cite{koornw}. For an inductive proof, see
Slater~\cite[Sec.~3.3.1]{slater}. For another elementary classical proof,
see Gasper and Rahman~\cite[Sec.~2.6]{grhyp}.

Concluding this section, we would like to add another thought,
kindly initiated by an anonymous referee.
It is worth comparing our proof with Askey and Ismail's~\cite{askmail}
elegant (and now classical) proof of Bailey's $_6\psi_6$ summation.
Their proof uses a method in this context often referred to as
``Ismail's argument"
since Ismail~\cite{ismail} was apparently the first to apply
Liouville's standard analytic continuation argument in the context of
bilateral basic hypergeometric series.
Askey and Ismail use Rogers' $_6\phi_5$ summation once to evaluate
the $_6\psi_6$ series at an infinite sequence and then apply
analytic continuation.
Here, we evaluate the $_6\psi_6$ series on a domain, 
and, for the full theorem, we also need analytic continuation.
In fact, we need, in addition to $|a^2q/bcde|<1$ two other inequalities
on $a,b,c,d,e$, namely $|aq/bc|<1$ and $|aq/cde|<1$, in order to apply 
the $_6\phi_5$ summation theorem. In the end, these additional conditions
can be removed. In particular, both sides of identity \eqref{66gl}
are analytic in $1/c$ around the origin. So far, we have shown
the identity for $|1/c|<\min(|b/aq|,|de/aq|,|bde/a^2q|)$.
By analytic continuation, we extend the identity to be valid
for $|1/c|<|bde/a^2q|$, the radius of convergence of the series.

\end{document}